\documentclass[11pt,a4paper]{article}

\usepackage[left=18mm,right=18mm,top=20mm,bottom=22mm]{geometry}
\usepackage[T1]{fontenc}
\usepackage{lmodern}
\usepackage{amsmath,amssymb,amsthm,mathtools}
\usepackage{xcolor}
\usepackage{enumitem}
\usepackage[expansion=false]{microtype}
\usepackage{hyperref}

\hypersetup{
  colorlinks=true,
  linkcolor=blue!45!black,
  citecolor=red,
  urlcolor=blue!45!black,
  pdftitle={Discreteness of critical points and fibres and unique continuation for conformal phi-curves and weakly-conformal k-harmonic maps},
  pdfauthor={Anton Iliashenko}
}

\setlist{leftmargin=*,itemsep=2pt,topsep=4pt}
\numberwithin{equation}{section}

\newtheorem{mainthm}{Theorem}

\newtheorem{theorem}{Theorem}[section]
\newtheorem{lemma}[theorem]{Lemma}
\newtheorem{proposition}[theorem]{Proposition}

\theoremstyle{definition}
\newtheorem{definition}[theorem]{Definition}
\newtheorem{example}[theorem]{Example}
\theoremstyle{remark}
\newtheorem{remark}[theorem]{Remark}

\newcommand{\R}{\mathbb{R}}
\DeclareMathOperator{\Ric}{Ric}
\DeclareMathOperator{\Crit}{Crit}
\DeclareMathOperator{\Inj}{Inj}
\DeclareMathOperator{\vol}{vol}
\newcommand{\II}{\mathrm{II}}
\newcommand{\fint}{\mathop{\mathchoice
 {\displaystyle\int\!\!\!\!-}{\textstyle\int\!\!\!-}
 {\scriptstyle\int\!\!-}{\scriptscriptstyle\int\!-}}}

\title{Discreteness of critical points and fibres\\
  and two-map unique continuation for conformal $\phi$-curves\\
  and weakly-conformal $k$-harmonic maps}
\author{Anton Iliashenko\\
  \small Beijing Institute of Mathematical Sciences and Applications\\
  \small Beijing, China\\
  \small\href{mailto:antoniliashenko@bimsa.cn}{\texttt{antoniliashenko@bimsa.cn}}}
\date{November 2026}

\begin{document}
\maketitle

\begin{abstract}
We study conformal $\phi$-curves and weakly-conformal
$k$-harmonic maps $u:L^k\to M^n$ between Riemannian manifolds,
where $L$ is connected, $3\le k\le n$ and $\phi$ is a semicalibration on $M$.
For nonconstant $C^1$ maps in either class, we prove that the
critical set and every fibre are discrete.
We show that smooth nonconstant maps are immersions.
We also prove that if two maps in the same class,
with respect to the same data, agree on a nonempty
open set, then they agree everywhere.
Finally, when $L$ is closed, we prove that
if the set of injective points is nonempty, then it is open and dense.
\end{abstract}

\tableofcontents

\section{Introduction}\label{sec:introduction}

Conformal $\phi$-curves, where $\phi$ is a (semi)calibration, combine weak conformality and a
calibration condition, which form a natural generalization of classical
$J$-holomorphic curves.
Away from the critical set, conformal $\phi$-curves parametrize
calibrated immersed submanifolds. The critical set is therefore the interesting
set which distinguishes this theory from the theory of just calibrated submanifolds.
On the other hand, the presence of critical points makes the analysis much harder.

A lot of progress has recently been done on this class of maps.
Cheng--Karigiannis--Madnick~\mbox{\cite{CKM}} prove various analytic properties,
analogous to those of $J$-holomorphic curves, such as interior regularity,
removability of singularities, energy gap, and bubble-tree convergence. They initially do it
just for those conformal $\phi$-curves for which $\phi$ comes from a vector cross product.
Iliashenko--Karigiannis~\mbox{\cite{IK}} observe that these results hold for any (semi)calibration.
For the associative and coassociative $G_2$-calibrations and the Cayley
$\mathrm{Spin}(7)$-calibration, Cheng--Karigiannis--Madnick~\cite{CKMVariational}
characterize calibrated submanifolds variationally and apply this
characterization to conformal $\phi$-curves.
Ikonen--Pankka~\cite{IP} address rigidity by studying which calibrations
admit a Liouville theorem. In the Euclidean setting, Ikonen~\cite{IkonenDiscrete}
proves discreteness of critical sets and fibres for nonconstant conformal
$\phi$-curves, assuming that the domain has dimension at least three and
the calibration has constant coefficients. There are also further generalizations, see
such as quasiregular curves, studied by Pankka~\cite{Pankka} and Ikonen--Pankka~\cite{IP}.

Recently, conformal $\phi$-curves were used by Broder--Iliashenko--Madnick~\cite{BIM} to
introduce notions of calibrated hyperbolicity analogous to those from classical complex geometry.
Instead of using holomorphic curves, they used conformal $\phi$-curves from $\R^k$ or the Poincar\'e
$k$-ball to define
generalizations of Kobayashi-Royden pseudo-metric, Brody hyperbolicity, Kobayashi hyperbolicity, and others.
Then, Iliashenko--Karigiannis--Madnick~\cite{IKM} prove calibrated analogues of
Royden's and Kiernan's theorems, relating hyperbolicity of the target to
equicontinuity and normality of families of conformal $\phi$-curves from
the Poincar\'e $k$-ball. They also obtain a corresponding Montel theorem.
Ikonen--Pim~\cite{IkonenPim} use a rescaling argument to prove
a calibrated analogue of Brody's theorem.

Weakly-conformal $k$-harmonic maps form a natural class in their own right. On their regular set,
these maps parametrize minimal submanifolds; see Takeuchi~\cite{Takeuchi}
and Chen--Wei~\cite{ChenWei}. The two classes are closely related: when $\phi$ is closed, every conformal
$\phi$-curve is $k$-harmonic~\cite{CKM,IK}. Without closedness, the two
classes still share the properties central to our arguments:
weak conformality and locally bounded mean curvature on the regular set.

We study critical points, fibres, and unique continuation for both classes,
allowing curved domains and targets. Assuming the domain is closed,
we also examine the set of injective points.

Let $(L^k,g)$ and $(M^n,h)$ be smooth
Riemannian manifolds without boundary, where $L$ is connected.
We assume $n\ge k\ge3$ and consider $C^1$ maps $u:L\to M$. When $L$ is oriented, $\vol_g$
denotes its volume form.
The two classes we consider are:
\begin{enumerate}[label=(\roman*)]
\item A \emph{conformal $\phi$-curve}: assume $L$ is oriented,
$\phi$ is a smooth \emph{semicalibration} on $M$, and for some function $\lambda:L\to[0,\infty)$,
\begin{equation}\label{eq:smith}
  u^*h=\lambda^2g,\qquad u^*\phi=\lambda^k\vol_g.
\end{equation}
\item A \emph{weakly-conformal $k$-harmonic map}: for some function $\lambda:L\to[0,\infty)$,
\[
  u^*h=\lambda^2g, \qquad \operatorname{div}_g\bigl(|du|_{g,h}^{k-2}du\bigr)=0
  \quad\text{weakly}.
\]
\end{enumerate}
The first equation in both classes means that the map is~\emph{weakly-conformal}
and $\lambda$ has to be $\frac{|du|_{g,h}}{\sqrt{k}}$.

The first group of theorems concerns the \emph{critical set}, defined by
\[
  \Crit(u)\coloneqq\{p\in L:du_p=0\}.
\]

\begin{mainthm}[Discreteness of the critical set]\label{thm:A}
For a nonconstant map in either class, $\Crit(u)$ is discrete in $L$.
In particular, every compact subset of $L$ contains only finitely many
critical points.
\end{mainthm}

\begin{mainthm}[Absence of critical points for smooth maps]\label{thm:B}
A smooth map in either class is either constant or an immersion.
\end{mainthm}

\begin{mainthm}[Two-map unique continuation]\label{thm:C}
Let $u,v:(L,g)\to(M,h)$ be $C^1$ maps which are either both
conformal $\phi$-curves for the same semicalibration $\phi$, or
both weakly-conformal and $k$-harmonic.
If they agree on a nonempty open subset of $L$, then $u=v$ on $L$.
\end{mainthm}

For $y\in M$, the \emph{fibre} over $y$ is
$u^{-1}(y)\coloneqq\{p\in L:u(p)=y\}$.
The next two results describe the fibres and the set of
\emph{injective points}, defined by
\[
  \Inj(u)\coloneqq
  \{p\in L:du_p\ne0,\ u^{-1}(u(p))=\{p\}\}.
\]

\begin{mainthm}[Discreteness of fibres]\label{thm:D}
For a nonconstant map in either class, every fibre $u^{-1}(y)$,
$y\in M$, is closed and discrete in $L$.
Consequently every compact subset of $L$ intersects each fibre in only
finitely many points. If $L$ is compact, every fibre is finite.
\end{mainthm}

\begin{mainthm}[Density of injective points]\label{thm:E}
Let $u$ be a nonconstant map in either class, and assume that $L$
is closed.
If $\Inj(u)\ne\varnothing$, then $\Inj(u)$ is open and dense in $L$.
Moreover, every $p\in\Inj(u)$ has open neighborhoods $U\subset L$
of $p$ and $V\subset M$ of $u(p)$ such that $u|_U$ is a smooth
embedding and $u^{-1}(V)\subset U$.
\end{mainthm}

For smooth nonconstant maps, Theorem~\ref{thm:B} makes the
nonzero-differential condition in $\Inj(u)$ automatic.

\begin{remark}[A broader class]\label{rem:broader-class}
Both classes~(i) and~(ii) above are contained in the broader class
of $C^1$ weakly-conformal maps whose restriction to
$L\setminus\Crit(u)$ is a smooth immersion and whose mean curvature
is uniformly bounded on $U\setminus\Crit(u)$ for every relatively
compact open set $U\Subset L$.
Using Martino's argument~\cite[Theorem~1.2(2)]{Martino}, one can show that maps in this
broader class belong to $C^{1,\alpha}_{\mathrm{loc}}$ for some $\alpha>0$.
Theorems~\ref{thm:A}, \ref{thm:B} and~\ref{thm:D}, as well as
continuation relative to a constant map, extend to this class by similar
arguments, which we omit.
Theorems~\ref{thm:C} and~\ref{thm:E} use additional properties of the
image available to us in classes~(i) and~(ii), so we restrict our focus to them.
\end{remark}

\begin{remark}[Domains with boundary]\label{rem:boundary}
When $L$ has a smooth boundary, Theorems~\ref{thm:A}--\ref{thm:D}
apply to its connected dense interior
$L^\circ\coloneqq L\setminus\partial L$.
Thus critical sets and fibres of nonconstant maps are discrete in
$L^\circ$, and smooth nonconstant maps are immersions there.
Equality in Theorem~\ref{thm:C} extends to $L$ by continuity for
maps continuous up to the boundary. The finiteness conclusions hold on compact subsets of $L^\circ$;
they do not say anything about boundary critical points or accumulation of
critical points or fibres at the boundary, even if $L$ is compact.
Theorem~\ref{thm:E} remains restricted to closed domains.
\end{remark}

\begin{remark}[The $k=2$ case]\label{rem:k-two}
The $k=2$ analogues of Theorems~\ref{thm:A}, \ref{thm:C},
\ref{thm:D} and~\ref{thm:E} hold for both classes~(i) and~(ii).
For $J$-holomorphic curves, all four conclusions are classical;
see McDuff--Salamon~\cite{McDuffSalamon}.

For the full classes~(i) and~(ii), the same argument used for $J$-holomorphic curves proves
discreteness of the critical set~\cite[Lemma~2.4.1]{McDuffSalamon}
\cite[Section~2.3]{Wendl}.
Two-map continuation follows from classical elliptic unique
continuation~\cite{Aronszajn}.
The proofs of Theorems~\ref{thm:D} and~\ref{thm:E} below also apply
when $k=2$.
For harmonic maps, Theorems~\ref{thm:A} and~\ref{thm:C} do not require
weak conformality~\cite[Section~2.3]{Wendl} and~\cite[Theorem~1]{Sampson}.
For an almost complex structure $J$, $J$-holomorphic curves are automatically
weakly-conformal for compatible metrics and harmonic for K\"ahler targets.

Theorem~\ref{thm:B} fails for $k=2$: the smooth, weakly-conformal,
holomorphic and harmonic map $z\mapsto z^m$ has a critical point at $0$ for every integer
$m \geq 2$.
\end{remark}

\begin{example}\label{ex:ejiri-fibres}
We briefly discuss the Ejiri cone construction~\cite[Theorem~1.4 and \S~6.3]{IkonenDiscrete},
which produces a conformal $\psi$-curve $F:\R^4\to\R^7$, for the standard coassociative calibration
$\psi$ on $\R^7$. It has one critical point, which is $0 \in \R^4$.
Moreover, for any $0\ne p\in\R^4$, we have $\#F^{-1}(F(p))=24$.
This makes $\Inj(F)=\varnothing$.
In particular, since $\R^4$ is not compact and $\Inj(F)=\varnothing$,
Theorem~\ref{thm:E} does not apply.
\end{example}

The author in~\mbox{\cite[Theorem~B]{Branch}} proves Taylor flatness at critical
points for smooth weakly-conformal maps. There, under an additional local curvature
bound on the scalar curvature of $u^{*}h$, he excludes critical points for
nonconstant maps. Theorem B of this paper
requires no such curvature assumption, but our maps are more than just weakly-conformal.
Hence, for general smooth weakly-conformal maps with $n>k$, the analogue of Theorem B remains open.
Recall that in equal dimensions, it is indeed proven by
Baird--Wood~\cite[Theorem~11.4.6]{BW}.

Ikonen~\mbox{\cite[Theorem~1.1]{IkonenDiscrete}} proves that the critical
set of a nonconstant $C^1$ conformal $\phi$-curve between Euclidean spaces
is discrete, where $\phi$ is a constant-coefficient calibration. Note that the
last condition implies that $\phi$ is closed, which makes the conformal
$\phi$-curve also $k$-harmonic. Thus,
Theorem~\mbox{\ref{thm:A}} extends his conclusion to any smooth Riemannian
domain and target manifold, to any smooth semicalibration, and separately to any
weakly-conformal $k$-harmonic map.
Theorem~\ref{thm:D} extends his discrete-fibre
theorem~\cite[Theorem~1.3]{IkonenDiscrete} to the same settings.
Pankka~\cite[Corollary~1.15]{Pankka} claims discreteness of fibres for nonconstant $C^1$
quasiregular curves,
but, as noted in personal communication, the argument for this corollary contains a gap.
In particular, the step that produces a single coordinate projection whose
composition with the quasiregular curve remains quasiregular on a
neighborhood of a critical point is not justified.
The same gap also affects the proof that $du\ne0$ almost everywhere and the claim
that every $C^1$ quasiregular curve is locally a graph over a quasiregular map.

Theorem~\ref{thm:C} extends Ikonen's continuation relative to a
constant~\cite[Theorem~1.2]{IkonenDiscrete} to comparison of two maps
in these settings. Note that Theorem~\ref{thm:C}
does not establish strong two-map continuation
at a shared critical point. The author is currently investigating this question.
The density conclusion in Theorem~\ref{thm:E} gives a
higher-dimensional analogue of the classical result for
$J$-holomorphic curves~\cite[Proposition~2.5.1]{McDuffSalamon}.

We now outline the proof strategy. Its starting point is
a classical Gauss-equation estimate (Lemma~\ref{lem:gauss}), which,
for the induced metric $\gamma\coloneqq u^*h$ on $L\setminus\Crit(u)$,
gives a Ricci bound $\Ric_\gamma\le C\gamma$ uniform near
critical points for the two classes considered here.
For a conformal $\phi$-curve, the regular mean curvature is locally
controlled by $d\phi$; for a weakly-conformal $k$-harmonic map it
vanishes.

For Theorem~\mbox{\ref{thm:A}}, we apply Ma--Qing's~\mbox{\cite{MQHuber}}
technique to construct a local Riesz measure~\mbox{\cite{HKM}} on a neighbourhood of a chosen
critical point.
This measure assigns uniformly positive masses to each critical point in the same neighborhood.
More precisely, after a local conformal change of the metric, we can make $-\log \lambda$
$k$-superharmonic.
This step is also classical~\mbox{\cite{MQPotential}}. After this,
we average an estimate of Kilpel\"ainen--Mal\'y~\mbox{\cite{KM94}} involving the Wolff potential.
Since the Riesz measure is finite on compact sets, this gives that the critical points are discrete.

For Theorem~\mbox{\ref{thm:B}}, we prove that the dilation of a nonconstant map cannot vanish to
infinite order at a critical point.
Together with the Taylor-flatness theorem from~\mbox{\cite{Branch}}, this excludes critical
points for smooth nonconstant maps.
For Theorem~\ref{thm:C}, we first compare regular image graphs
using classical elliptic unique continuation and then identify
their conformal parametrizations by classical two-jet rigidity. For Theorem~\mbox{\ref{thm:D}},
we apply the monotonicity inequality for smooth
immersions with boundary~\mbox{\cite{Scharrer}} on
annuli around a critical point. The resulting area bound rules out
accumulation of regular preimages there.
For Theorem~\mbox{\ref{thm:E}}, we remove the fibres over critical values and
apply unique continuation to the graphs of $u$ near the remaining fibres.

Our arguments differ substantially from Ikonen's.
He obtains a conservation law from the natural differential,
constructs a monotone frequency, and analyzes homogeneous blow-ups.

\subsection*{Acknowledgements}
The author acknowledges the support of the International Scientists Project BJNSF \#25ISA023.

\subsection*{AI disclosure}
The author used AI (Artificial Intelligence) to help with literature searches and to
develop and check the proofs. The manuscript is written by the author.

\subsection*{Organization}
Section~\mbox{\ref{sec:notation}} fixes notation and conventions and recalls needed standard
results.
Section~\mbox{\ref{sec:preliminaries}} recalls regularity of our maps and establishes the
required curvature estimates.
Section~\mbox{\ref{sec:critical-set}} proves that the critical set of a nonconstant
map is discrete (Theorem~\mbox{\ref{thm:A}}). Section~\mbox{\ref{sec:smooth}} shows that
every smooth nonconstant map in either class has no critical points (Theorem~\mbox{\ref{thm:B}}).
Section~\mbox{\ref{sec:continuation}} proves unique continuation for two maps
(Theorem~\mbox{\ref{thm:C}}). Section~\mbox{\ref{sec:fibres}} proves that every fibre of a
nonconstant map is discrete (Theorem~\mbox{\ref{thm:D}}). Section~\mbox{\ref{sec:injective}}
proves that when the domain is closed, if the set of injective points is nonempty, then it is
open and dense (Theorem~\mbox{\ref{thm:E}}).

\section{Notation and background results}\label{sec:notation}

We fix notation and conventions.

\subsection{Preliminaries}
Let $(L^k,g)$ and $(M^n,h)$ be smooth Riemannian manifolds
without boundary, where $L$ is connected and $n\ge k\ge3$. A \emph{closed manifold} means
compact with no boundary.
We write $u:L\to M$ for a $C^1$ map and call $L$ its \emph{domain}
and $M$ its \emph{target}.
We write $\operatorname{dist}_g$ for Riemannian distance on $L$ and
$\operatorname{div}_g$ for \emph{divergence}.

For $U \subset L$, we write $U\Subset L$ to mean that $U$ is open in $L$
and has compact closure in $L$.

We write $B_g(p,r)\coloneqq\{q\in L:\operatorname{dist}_g(p,q)<r\}$ and
$B(x,r):=\{y\in\mathbb R^k:|x-y|<r\}$.

For a measurable set $E\subset\R^k$, $|E|$ denotes its
\emph{Lebesgue measure}.

Let $f:E\to\R^N$, where $N\ge1$ is an integer.
If $0<|E|<\infty$ and $f$ is integrable on $E$, we define its
\emph{average} by
\[
  \fint_E f\,dx \coloneqq\frac1{|E|}\int_E f\,dx.
\]

A set $E\subset L$ is \emph{discrete} if each of its points is isolated.
$\#F$ denotes the number of elements of a finite set $F$.
For $y\in M$, the \emph{fibre} is
$u^{-1}(y)\coloneqq\{p\in L:u(p)=y\}$.

For domains with boundary as in Remark~\mbox{\ref{rem:boundary}},
the \emph{interior} $L^\circ \coloneqq L \setminus \partial L$ of $L$ is connected and dense in $L$.

\subsection{Function spaces}
Fix an integer $j \geq 0$. When the domain is understood, $C^j$ denotes a set of maps
whose coordinate derivatives of order $0$ to $j$ exist and are continuous.
A \emph{multi-index} is a vector
$\beta\in\{0,1,2,\ldots\}^k$, whose \emph{total order} is
$|\beta|\coloneqq\sum_{i=1}^k\beta_i$, and we write
$\partial^\beta\coloneqq\partial_1^{\beta_1}\cdots\partial_k^{\beta_k}$.
We use the Einstein summation convention to sum over repeated indices.

Let $U\subset\R^k$ be open and let $f:U\to\R^N$, where $N \geq 1$ is an integer.
For a compact set $K\subset U$ and $0<\alpha\le1$, set
\[
  [f]_{\alpha;K}\coloneqq
  \sup_{\substack{x,y\in K\\x\ne y}}
  \frac{|f(x)-f(y)|}{|x-y|^\alpha}.
\]

When this quantity is finite, we say that $f$ is
\emph{H\"older continuous} with exponent $\alpha$ on $K$. We say that $f$ is
\emph{locally H\"older continuous} with exponent
$\alpha$ on $U$ if $[f]_{\alpha;K}<\infty$ for every compact
$K\subset U$. The \emph{H\"older space} $C^{j,\alpha}_{\mathrm{loc}}(U)$ consists of $C^j$ functions
whose order-$j$ derivatives are locally H\"older continuous with exponent $\alpha$.
For a map $u$ between manifolds, we ask for the same conditions in smooth
local coordinates. Clearly these notions are independent of the choice of coordinates.
We write $du$ for the differential. Now, fix local coordinates. For differentiable $u$, we
define its \emph{Jacobian matrix} by
$Du\coloneqq(\partial_i u^\alpha)_{\alpha,i}$.
For an integer $b\ge2$ and $u \in C^b$,
$D^b u \coloneqq \bigl(\partial^\beta u^\alpha\bigr)_{\alpha,\ |\beta|=b}$.
The \emph{two-jet} of a $C^2$ map $u$ at a point $p$ consists of $u(p)$, $du(p)$ and $D^2 u(p)$.

For $1\le q<\infty$, $L^q$ denotes the space of measurable functions
$f:L \rightarrow \mathbb{R}$ with $\int_L |f|^q d\vol_g < \infty.$
For an open set $U\subset\R^k$, let $C_c^\infty(U)$ denote the smooth
functions with compact support in $U$.
The \emph{Sobolev space} $W^{1,k}(U)$ consists of functions $f\in L^k(U)$
for which there are $F_1,\ldots,F_k\in L^k(U)$ satisfying
\[
  \int_U f\,\partial_i\eta\,dx=-\int_U F_i\eta\,dx
  \qquad\forall\eta\in C_c^\infty(U),\quad i=1,\ldots,k.
\]
For $f\in W^{1,k}(U)$, we define its \emph{weak derivatives} by
$\partial_i f\coloneqq F_i$. A subscript $\mathrm{loc}$ means that the stated
condition holds on every compact subset of the domain. Vector-valued spaces
are defined componentwise.
We use elements of $W^{1,k}_0(U) \coloneqq\overline{C_c^\infty(U)}^{\,W^{1,k}(U)}$ as
\emph{test functions}.

Let $A\in L^{\frac{k}{k-1}}_{\mathrm{loc}}(U;\R^k)$ and
$F\in L^{\frac{k}{k-1}}_{\mathrm{loc}}(U)$. We say that
$\sum_{i=1}^k\partial_iA_i=F$ holds \emph{weakly} if
\[
  -\int_\Omega\sum_{i=1}^k A_i\,\partial_i\eta\,dx
  =\int_\Omega F\eta\,dx
  \qquad\forall\Omega\Subset U,\quad\forall\eta\in W^{1,k}_0(\Omega).
\]

The letter $C$ in inequalities denotes a constant that might change from line to line.

Fix $p \in L$ and set $r(x)\coloneqq\operatorname{dist}_g(x,p)$ and let
$f:L \rightarrow \mathbb R^N$.
For $a > 0$, we write $f=O(r^a)$ to mean that there exist $C,\varepsilon > 0$
such that $|f(x)|\le Cr^a$ if $0 < r(x) < \varepsilon$. We say $f$
\emph{vanishes to infinite order} at $p$ if
$f = O(r^a)$ for every $a > 0.$ If $f:L\rightarrow M$ is smooth, we say that $f$ is
\emph{Taylor flat} at $p$ if in local coordinates, $D^b f (p) = 0$ for all $b = 1,2,3,\dots.$
If $f_1,f_2:L \rightarrow M$ are two $C^1$ maps, then we say $f_1$ and $f_2$
\emph{agree to infinite order} at $p$ if
in local coordinates, $|f_1-f_2| = O(r^a)$ for all $a > 0.$ If $f_1,f_2$ are smooth,
the last condition is equivalent to $f_1(p)=f_2(p)$ and $D^b f_1(p)=D^b f_2(p)$ for every
integer $b\ge1$, in local coordinates.

\subsection{Weak conformality and curvature}\label{sec:weak-conformality}
The norm $|du|_{g,h}$ is the \emph{Hilbert--Schmidt norm}:
for a $g$-orthonormal frame $e_1,\ldots,e_k$,
\[
  |du|_{g,h}^2\coloneqq\sum_{i=1}^k h(du(e_i),du(e_i)),\qquad
  (u^*h)(X,Y)\coloneqq h(du(X),du(Y)).
\]
An \emph{immersion} has injective differential at every point.
A \emph{smooth embedding} is an injective immersion that is a homeomorphism
onto its image.

For the maps in this paper, we say $p \in L$ is a \emph{critical point} if $du_p = 0$; otherwise
we call it a~\emph{regular point}.
The \emph{critical set}
and the \emph{regular set} are defined correspondingly as
\[
  Z_u:=\Crit(u):=\{p\in L:du_p=0\},
  \qquad
  \mathcal R_u:=L\setminus Z_u.
\]

We call $p\in L$ an \emph{injective point} of $u$ if $du_p\ne0$
and $u^{-1}(u(p))=\{p\}$, following the terminology of
McDuff--Salamon~\cite[Section~2.5]{McDuffSalamon}.
We denote the set of such points by $\Inj(u)$.

Set $\lambda\coloneqq\frac{|du|_{g,h}}{\sqrt{k}}\ge0$.
A map is \emph{weakly-conformal} if $u^*h=\lambda^2g$.
Its differential then has rank $k$ or is zero, and consequently
$Z_u = \{\lambda=0\}$ and
$\mathcal R_u = L\setminus \Crit(u) = \{\lambda>0\}$.

The \emph{induced metric} $\gamma\coloneqq u^*h=\lambda^2g$ is positive definite
on $\mathcal R_u$. Geometric quantities involving $\gamma$ are used
only on this set.
For a metric $q$, the coefficients $q^{ij}$ are the entries of the
inverse matrix of $(q_{ij})$.
We use $\vol_g$ for the oriented volume form and $d\vol_g$ for
the Riemannian volume measure, which also exists without an orientation.

Let $V\subset\mathcal R_u$ be an open set on which $u$ is a smooth
embedding, and set $\Sigma\coloneqq u(V)\subset M$.
Write $\nabla^h$ for the Levi--Civita connection of $h$.
The \emph{normal bundle} is
$N\Sigma\coloneqq(T\Sigma)^{\perp_h}\subset TM|_\Sigma$,
and $\perp$ denotes orthogonal projection onto it.
We identify $\gamma|_V$ with the induced metric on $\Sigma$ via $u|_V$.
For a vector bundle $E$, let $\Gamma(E)$ denote its space of smooth
sections.
For $A,B\in\Gamma(T\Sigma)$, the \emph{second fundamental form}
and \emph{averaged mean curvature} are defined by
\[
  \begin{aligned}
    \II(A,B)&\coloneqq(\nabla^h_A B)^\perp\in\Gamma(N\Sigma),\\
    H&\coloneqq\frac1k\operatorname{tr}_\gamma\II\in\Gamma(N\Sigma).
  \end{aligned}
\]
We also use $\II$ for the pullback second fundamental form on $V$:
for $p\in V$ and $v,w\in T_pL$, set
\[
  \II_p(v,w)\coloneqq\II_{u(p)}(du_p(v),du_p(w)).
\]
For $y\in\Sigma$ and a positively oriented $h$-orthonormal basis
$E_1,\ldots,E_k\in T_y\Sigma$, set
$\tau(y)\coloneqq E_1\wedge\cdots\wedge E_k$; this is the oriented
unit tangent $k$-vector at $y$.
A \emph{local bound for $H$} means a single bound on $U\cap\mathcal R_u$
for each relatively compact domain neighborhood $U$.
For a metric $q$ with Levi--Civita connection $\nabla^q$ and smooth
vector fields $X,Y,Z,W$, let $[X,Y]$ denote the Lie bracket of $X$ and $Y$.
We define the \emph{curvature operator} and \emph{Riemann curvature tensor} by
\[
  \begin{aligned}
    R^q(X,Y)Z&\coloneqq
    \nabla^q_X\nabla^q_YZ-\nabla^q_Y\nabla^q_XZ-\nabla^q_{[X,Y]}Z,\\
    \operatorname{Rm}_q(X,Y,Z,W)&\coloneqq q(R^q(X,Y)Z,W).
  \end{aligned}
\]
In a local frame $e_i$, set
$\operatorname{Rm}_{ijk\ell}(q)\coloneqq
\operatorname{Rm}_q(e_i,e_j,e_k,e_\ell)$.
We define the \emph{Ricci tensor} by contraction:
\[
  (\Ric_q)_{jk}\coloneqq q^{i\ell}\operatorname{Rm}_{ijk\ell}(q).
\]
For $y\in M$ and a two-plane $E\subset T_yM$, we write $\sec_h(E)$
for its sectional curvature with respect to $h$.
For linearly independent $A,B\in T_yM$, we also write
\[
  \sec_h(A\wedge B)\coloneqq\sec_h(\operatorname{span}\{A,B\}).
\]
If $A,B$ are $h$-orthonormal, then
$\sec_h(A\wedge B)=\operatorname{Rm}_h(A,B,B,A)$.

For a constant $C>0$, when we write $\Ric_\gamma\le C\gamma$ we mean that
\[
  (\Ric_\gamma)_p(V,V)\le C\gamma_p(V,V)
  \qquad(p\in\mathcal R_u,\ V\in T_pL).
\]

\subsection{The \texorpdfstring{$k$}{k}-Laplacian and the Riesz measure}
Let $U\subset L$ be open and $v:U\to\R$ be smooth.
We write $\nabla_gv$ for the \emph{gradient} of $v$,
$\nabla_g^2v$ for its \emph{Hessian}, and
$\Delta_g v\coloneqq\operatorname{div}_g(\nabla_g v)$ for its \emph{Laplacian}. The
\emph{$k$-Laplacian} of $v$ is defined as
\[
  \Delta_{k,g}v
  \coloneqq\operatorname{div}_g\bigl(|\nabla_gv|^{k-2}\nabla_gv\bigr).
\]
For a scalar function $v\in W^{1,k}_{\mathrm{loc}}(U)$, its \emph{weak gradient} is
$\nabla_gv\coloneqq\sum_{i,j=1}^k g^{ij}(\partial_jv)\partial_i$,
where $\partial_jv$ are its weak derivatives.
We call $v$ a \emph{$k$-harmonic function} if,
for every open $\Omega\Subset U$,
\[
  \int_\Omega |\nabla_gv|_g^{k-2}
  \langle\nabla_gv,\nabla_g\eta\rangle_g\,d\vol_g=0
  \qquad\forall\eta\in W^{1,k}_0(\Omega).
\]
A function $v\in W^{1,k}_{\mathrm{loc}}(U)$ is a \emph{weak supersolution} if,
for every open $\Omega\Subset U$,
\[
  \int_\Omega |\nabla_gv|_g^{k-2}
  \langle\nabla_gv,\nabla_g\eta\rangle_g\,d\vol_g\ge0
  \qquad\forall\eta\in W^{1,k}_0(\Omega),\ \eta\ge0.
\]
For smooth $v$, this condition is equivalent to $\Delta_{k,g}v\le0$
pointwise, and in this case we call $v$ a \emph{smooth supersolution}.

For a $C^1$ map $u:U\to M$, we use the pullback connection. For smooth $u$, we set
\[
  \tau_k(u)\coloneqq\operatorname{div}_g(|du|^{k-2}du).
\]
For $u:U\to M$, let $C_c^1(U,u^*TM)$ denote the
compactly supported $C^1$ fields $X$ with $X(p)\in T_{u(p)}M$.
We say that $u$ is \emph{weakly $k$-harmonic} (or just \emph{$k$-harmonic}), and write
$\tau_k(u)=0$ \emph{weakly}, if
\[
  \int_U |du|_{g,h}^{k-2}
  \langle du,\nabla X\rangle_{g,h}\,d\vol_g=0
  \qquad\forall X\in C_c^1(U,u^*TM).
\]
Equivalently, for every open $\Omega\Subset U$, $u$ is a
\emph{stationary point} of the \emph{$k$-energy}
\[
  E_k(u;\Omega)\coloneqq\frac{1}{\sqrt{k}^k}\int_\Omega |du|_{g,h}^k\,d\vol_g
\]
under compactly supported variations of $u$ in $\Omega$.

A function $v:U\to(-\infty,+\infty]$ is
\emph{lower semicontinuous} if
$v(x)\le\liminf_{i\to\infty}v(x_i)$ whenever $x_i\to x$.
Such a function is \emph{$k$-superharmonic} if it is finite at some point
of each connected component of $U$ and, for every open set $V\Subset U$
and every continuous function $w:\overline V\to\R$ that is
$k$-harmonic on $V$, $w\le v$ on $\partial V$ implies $w\le v$ on $V$.

A \emph{Radon measure} is a nonnegative Borel measure that is finite on
compact sets and both inner and outer regular. For a Radon measure $\mu$ on $U$
and a Borel set $E\subset U$, we call $\mu(E)$ the \emph{mass} of $E$.
We say that $\mu$ has an \emph{atom} at $p\in U$ if $\mu(\{p\})>0$.

For a relatively compact open set $D\subset L$ and a compact set
$E\subset D$, the \emph{relative $k$-capacity} is
\[
  \operatorname{cap}_{k,g}(E,D)
  \coloneqq\inf_{\substack{\eta\in C_c^\infty(D)\\
  \eta\ge1\text{ near }E}}
  \int_D|\nabla_g\eta|^k\,d\vol_g,
\]
where by ``$\eta\ge1$ near $E$'' we mean that there exists an open set
$O$ with $E\subset O\subset D$ such that $\eta(x)\ge1$ for every $x\in O$.
For an open set $O\subseteq D$ and an arbitrary set $A\subseteq D$, we set
\[
  \begin{aligned}
    \operatorname{cap}_{k,g}(O,D)
    &\coloneqq\sup_{\substack{K\subseteq O\\K\text{ compact}}}
    \operatorname{cap}_{k,g}(K,D),\\
    \operatorname{cap}_{k,g}(A,D)
    &\coloneqq\inf_{\substack{A\subseteq O\subseteq D\\O\text{ open}}}
    \operatorname{cap}_{k,g}(O,D).
  \end{aligned}
\]
A set $E\subset L$ has \emph{locally zero $k$-capacity} if, for every pair of
open sets $D'\Subset D\Subset L$,
\[
  \operatorname{cap}_{k,g}(E\cap D',D)=0.
\]
Two smooth Riemannian metrics $g_1,g_2$ on $D$ are
\emph{uniformly equivalent} if there is a constant $C\ge1$ such that
\[
  C^{-1}(g_1)_p(X,X)\le(g_2)_p(X,X)\le C(g_1)_p(X,X)
  \qquad(p\in D,\ X\in T_pL).
\]
Such metrics give the same capacity-zero sets in $D$.

The following standard result supplies the Riesz measure and the estimate
used to prove discreteness of the critical set.

\begin{theorem}\label{thm:potential}
Fix a connected open set $D\subset\R^k$ with a metric $g$, and 
a $k$-superharmonic function $v:D \rightarrow [0,\infty]$. Then $v$
is finite almost everywhere and there exists a Radon measure
$\mu\coloneqq-\Delta_{k,g}v$, which we call a \emph{Riesz measure}.
Moreover, $\{v=+\infty\}$ has locally zero $k$-capacity.
Finally, for any $x\in D$ and $R>0$ such that $B(x,3R)\Subset D$, we have
\begin{equation}\label{eq:wolff}
  v(x)\le C\inf_{B(x,R)}v
  +C\int_0^{2R}\mu(B(x,t))^{\frac{1}{k-1}}\,\frac{dt}{t},
\end{equation}
for some constant $C > 0$ which depends only on $k$ and $g|_{B(x,3R)}$.
The integral in \mbox{\eqref{eq:wolff}} is called the \emph{Wolff potential}.
\end{theorem}
\begin{proof}
In local coordinates, the metric $k$-Laplacian has divergence form
\[
  \Delta_{k,g}v
  =\frac1{\sqrt{\det g}}\,
  \partial_i\!\left(\sqrt{\det g}\,|dv|_g^{k-2}g^{ij}\partial_jv\right).
\]
On each relatively compact subset, smoothness and positive definiteness
of $g$ give the conditions for this operator to be in the class considered by
Heinonen--Kilpel\"ainen--Martio~\cite[(3.3)--(3.7)]{HKM}.
Their Theorem~21.2 gives the Riesz measure.
The implication (i)$\Rightarrow$(iii) in their Theorem~10.1 gives locally zero capacity of
$\{v=+\infty\}$;
capacity-zero sets have Lebesgue measure zero, so $v$ is finite almost everywhere.
Finally, we obtain \eqref{eq:wolff} from
Kilpel\"ainen--Mal\'y~\cite[Theorem~1.6]{KM94} with both exponent and dimension equal to $k$.
\end{proof}

\subsection{Semicalibrations and conformal \texorpdfstring{$\phi$}{phi}-curves}
The following definitions specify the additional condition in class~(i).
The space $\Omega^k(M)$ consists of smooth differential $k$-forms.
The \emph{comass} of $\phi \in \Omega^k(M)$ at $y\in M$ is
\[
  \|\phi_y\|_h^*
  \coloneqq\sup\{|\phi_y(v_1,\ldots,v_k)|:
  v_1,\ldots,v_k\text{ are }h\text{-orthonormal}\}.
\]
A \emph{semicalibration} has comass at most one; a \emph{calibration} is additionally
\emph{closed}, meaning $d\phi=0$, where $d$ is the exterior derivative.
An oriented $k$-plane is \emph{$\phi$-calibrated} if $\phi$ takes
the value one on an oriented $h$-orthonormal basis of this plane.
The \emph{pullback} is defined by
$u^*\phi(X_1,\ldots,X_k)\coloneqq\phi(du(X_1),\ldots,du(X_k))$,
where $X_1,\ldots,X_k$ are smooth vector fields on $L$.

\begin{definition}\label{def:phi-curve}
Let $L$ be oriented and let $\phi\in\Omega^k(M)$ be a smooth
semicalibration. A $C^1$ map $u:(L,g)\to (M,h)$ is a \emph{conformal
$\phi$-curve} if, for some function $\lambda:L\to[0,\infty)$,
\[
  u^*h=\lambda^2g,\qquad u^*\phi=\lambda^k\vol_g.
\]
The first equation implies $\lambda=\frac{|du|_{g,h}}{\sqrt{k}}$.
Conversely, the second equation with $\lambda=\frac{|du|_{g,h}}{\sqrt{k}}$ implies weak
conformality of $u$ by the comass and
arithmetic--geometric mean inequalities~\mbox{\cite{CKM,IK}}.
At the points where $du_p \neq 0$, $du_p(T_pL)$, with orientation induced by $L$,
is $\phi$-calibrated.
\end{definition}

\begin{remark}\label{rem:smith-terminology}
We use the terminology of Ikonen--Pankka~\mbox{\cite{IP}}.
This class of maps is known under other names: \emph{Smith map}~\mbox{\cite{CKM}} and
\emph{multiholomorphic map}~\mbox{\cite{Smith}} are used when the semicalibration comes from a
vector cross product.
For general calibrations, the name
\emph{Smith immersion} is used in~\mbox{\cite{IK}}.
Note that since these are just different names for the same object, \emph{Smith immersion} does not
mean it is an immersion everywhere.

If $d\phi=0$, every conformal $\phi$-curve is $k$-harmonic~
\mbox{\cite{CKM,IK}}.
\end{remark}

\section{Regularity and curvature estimates}\label{sec:preliminaries}

First, we quickly recall the regularity of our maps. Then we obtain a bound on their mean curvature.

\begin{lemma}\label{lem:regularity}
A map in either class~(i) or~(ii)  belongs to
$C^{1,\alpha}_{\mathrm{loc}}$ for some $\alpha>0$
and is smooth on $\mathcal R_u$. Let $V\Subset\mathcal R_u$ be such that
$u|_V$ is an embedding, and set $\Sigma\coloneqq u(V)$.
Then $H$ satisfies
\begin{itemize}
  \item in class~(i), $-k\langle H,X\rangle=d\phi(X,\tau)$, for every $X \in N\Sigma.$
  \item in class~(ii), $H = 0$.
\end{itemize}
In particular, $|H|$ is uniformly bounded on $U\cap\mathcal R_u$
for every relatively compact open set $U\Subset L$.
\end{lemma}

\begin{proof}

For class (i), Cheng--Karigiannis--Madnick~\mbox{\cite[Section~1.2 and Theorem~4.10]{CKM}}
prove $C^{1,\alpha}_{\mathrm{loc}}$ regularity and smoothness on $\mathcal{R}_u$ for
semicalibrations coming from a vector cross product.
Iliashenko--Karigiannis~\mbox{\cite{IK}} note
that this result extends to any semicalibration.

For class~(ii), $C^{1,\alpha}_{\mathrm{loc}}$ regularity follows
from standard theory for $C^1$ $k$-harmonic
maps~\mbox{\cite{Uhlenbeck,Martino}}. Smoothness on $\mathcal R_u$ follows
from uniform ellipticity away from $du=0$.

For class~(i), take a compactly supported normal variation vector field $X$ on $\Sigma$,
and denote the varied images by $\Sigma_t$.
The quantity
$\operatorname{Vol}(\Sigma_t)-\int_{\Sigma_t}\phi$ has a minimum at $t=0$.
Taking the first variation yields
\[
  -k\int_\Sigma h(H,X)\,d\vol_\gamma
  =\int_\Sigma d\phi(X,\tau)\,d\vol_\gamma,
\]
where $\tau$ is the oriented unit tangent $k$-vector of $\Sigma$, as in
Section~\ref{sec:weak-conformality}.
Arbitrariness of $X$ gives
$-kh(H,X)=d\phi(X,\tau)$, and consequently
$|H|\le \frac{|d\phi|_h}{k}$.

For class~(ii), $u|_{\mathcal R_u}$ is a smooth conformal
$k$-harmonic immersion. Takeuchi~\cite[Corollary~4]{Takeuchi}
therefore gives $H=0$.

For every $U\Subset L$, continuity of $u$
makes $u(\overline U)$ compact. In class~(i), smoothness of $d\phi$
therefore gives a uniform bound for $|d\phi|_h$ on $u(\overline U)$ and hence
for $|H|$ on $U\cap\mathcal R_u$; in class~(ii), $H=0$.
\end{proof}

An upper bound for the ambient sectional curvature and a bound for $|H|$ give
an upper bound for $\Ric_\gamma$ by the Gauss equation. In the special case
$\sec_h \equiv \mathrm{Const}$ this is the Chen--Ricci
inequality~\mbox{\cite[Theorem~4]{ChenRicci}}.
We recall the argument.

\begin{lemma}\label{lem:gauss}
Let $u$ belong to either class~(i) or~(ii) and let $U\Subset L$ be open.
Choose $K\in\R$ and $H_0\ge0$ such that
\[
  \begin{aligned}
    \sec_h(E)&\le K,
    &&\text{for $y\in u(U),\ E\subset T_yM,\ \dim E=2,$}\\
    |H(p)|_h&\le H_0,
    &&\text{for $p\in U\cap\mathcal R_u.$}
  \end{aligned}
\]
Then, on $U\cap\mathcal R_u$, we have
\begin{equation}\label{eq:gauss}
  \Ric_\gamma\le
  \left((k-1)K+\frac{k^2H_0^2}{4}\right)\gamma.
\end{equation}
\end{lemma}
\begin{proof}
Since $u(\overline U)$ is compact, the ambient sectional curvature
has a finite upper bound there. Lemma~\ref{lem:regularity} provides $H_0$.
Fix a point $p\in U\cap\mathcal R_u$ and choose a
$\gamma$-orthonormal frame $e_1,\ldots,e_k$ of $T_pL$.
Set $e\coloneqq e_1$ and $E_i\coloneqq du_p(e_i)\in T_{u(p)}M$;
$E_1,\ldots,E_k$ form an $h$-orthonormal frame of
$du_p(T_pL)$. The contracted Gauss equation~\cite[Theorem~8.5]{LeeRM}
gives
\begin{equation}\label{eq:gauss-contracted}
  \begin{aligned}
    \Ric_\gamma(e,e)
    &=\sum_{i=2}^k\sec_h(E_1\wedge E_i)
    +k h(H,\II(e,e))-\sum_{i=1}^k|\II(e,e_i)|^2\\
    &\le (k-1)K+\Bigl(k h(H,\II(e,e))-|\II(e,e)|^2\Bigr)
    -\sum_{i=2}^k|\II(e,e_i)|^2.
  \end{aligned}
\end{equation}
Separately, we have
\begin{equation}\label{eq:gauss-square}
  k h(H,\II(e,e))-|\II(e,e)|^2
  =\frac{k^2}{4}|H|^2-
  |\II(e,e)-\tfrac k2H|^2
  \le\frac{k^2H_0^2}{4}.
\end{equation}
Substituting \eqref{eq:gauss-square} into \eqref{eq:gauss-contracted}
proves~\eqref{eq:gauss}.
\end{proof}

We state the intermediate results in a general form for a nonnegative function $\lambda$ on $(L,g)$,
which we later apply to the dilation of a weakly-conformal map. First, we need to
do a local conformal change of the metric to produce a Ricci gap.
The proof is elementary and uses the
standard conformal Ricci formula.

\begin{lemma}\label{lem:background}
Fix $(L^k,g)$ with $k \geq 3$ and
let $\lambda \geq 0$ be continuous, smooth on the nonempty set $\mathcal R\coloneqq\{\lambda>0\}$.
Set $\gamma\coloneqq\lambda^2g$ on $\mathcal R$.
Fix $p\in L$ with $\lambda(p)=0$ and an open neighborhood $U\subset L$ of $p$.
Suppose there is a constant $C\ge0$ such that
$\Ric_\gamma\le C\gamma$ on $U\cap\mathcal R$.
Then there are an open neighborhood $V\Subset U$ of $p$, a smooth
function $\psi:V\to\R$, and a constant $c_0>0$ such that $\widetilde g\coloneqq e^{2\psi}g$ and
$\widetilde\lambda\coloneqq e^{-\psi}\lambda$ satisfy
$\gamma=\widetilde\lambda^2\widetilde g=\lambda^2g$ and
\begin{equation}\label{eq:ricci-gap}
  Q\coloneqq\Ric_{\widetilde g}-\Ric_\gamma\ge c_0\widetilde g>0
  \quad\text{on }V\cap\mathcal R.
\end{equation}
On $V\cap\mathcal R$, the function $f\coloneqq\log\widetilde\lambda$ satisfies
\begin{equation}\label{eq:conformal-ricci}
  Q=(k-2)(\nabla_{\widetilde g}^2f-df\otimes df)
  +(\Delta_{\widetilde g}f+(k-2)|df|_{\widetilde g}^2)\widetilde g.
\end{equation}
\end{lemma}
\begin{proof}
Choose a normal neighborhood
$W\subset U$ of $p$ and define a smooth function
$\psi\coloneqq-a\,\operatorname{dist}_g(p,\cdot)^2$ on $W$, with $a>0$.
The conformal Ricci formula~\cite[Theorem~1.159(d)]{Besse} (note that Besse uses the opposite sign
for the Laplacian) at $p$ gives
\[
  \Ric_{\widetilde g}(p)
  =\Ric_g(p)+4(k-1)a\,g(p).
\]
Choose $a$ so that this is strictly positive definite. By continuity,
choose $c_0>0$ and an open neighborhood $V\Subset W$ of $p$ such that
$\Ric_{\widetilde g}\ge2c_0\widetilde g$ on $V$.
Since $\widetilde\lambda(p)=0$, shrink $V$ so that
$C\widetilde\lambda^2\le c_0$ on $V$. These inequalities prove
\eqref{eq:ricci-gap}.

The same conformal Ricci formula applied to $\gamma=e^{2f}\widetilde g$ on $V\cap\mathcal R$ gives
\eqref{eq:conformal-ricci}.
\end{proof}

\begin{remark}\label{rem:local-conformal-change}
We make this conformal change only on a neighborhood of a fixed
zero of $\lambda$. Within that local argument, we drop the tildes and,
without loss of generality, continue to denote the modified metric
and conformal factor by $g$ and $\lambda$.
\end{remark}

\section{Discreteness of the critical set}
\label{sec:critical-set}

We construct a nonnegative Radon measure on an
open neighborhood of a chosen critical point. We then show that this measure
assigns a uniformly positive mass to every critical point in a smaller
open neighborhood, thereby proving that the chosen critical point is
isolated.

The relation between logarithmic growth of a superharmonic function and
the mass of its Riesz measure appears in Ma--Qing~\cite{MQHuber}. In our
setting, we derive the needed bound by averaging the estimate in
Theorem~\ref{thm:potential}.

\begin{lemma}\label{lem:atom}
Let $g$ be a smooth metric on a connected open set $D\subset\R^k$,
and let $v\ge0$ be $k$-superharmonic on $D$.
Let $\mu\coloneqq-\Delta_{k,g}v$ denote its Riesz measure, as in
Theorem~\ref{thm:potential}. Fix $p\in D$. Suppose $\alpha>0$ and
there are $r_p>0$ and $C_p\in\R$ with $B(p,r_p)\subset D$ such that
\begin{equation}\label{eq:atom-log-lower}
  v(x)\ge\alpha\log\frac1{|x-p|}-C_p
  \qquad\text{for every }x\in B(p,r_p)\setminus\{p\}.
\end{equation}
Then $\mu(\{p\})\ge c\alpha^{k-1}$, where $c>0$ is independent
of $v$, $\alpha$, $r_p$ and $C_p$. We may choose $c>0$ uniformly for
$p$ in any fixed compact subset of $D$.
\end{lemma}

\begin{proof}
Fix $p_0\in D$ and choose $\rho>0$ such that
$\overline{B(p_0,8\rho)}\subset D$.
We prove the estimate uniformly for $p\in B(p_0,\frac{\rho}{4})$.
The change of coordinates $y\coloneqq\frac{x-p_0}{\rho}$ identifies $\overline{B(p_0,8\rho)}$
with $\overline{B(0,8)}$. We use the same symbols $D,g,v,\mu,p$ in these new coordinates and
allow $C_p$ to change.
We replace $r_p$ by $\frac{r_p}{\rho}$.
In these coordinates, $p\in B(0,\frac{1}{4})$.
By Theorem~\mbox{\ref{thm:potential}}, choose
$q\in B(0,\frac{1}{4})$ with $v(q)<\infty$. For $0<r<\frac{1}{8}$ set
\[
  A_r(p)\coloneqq B(p,2r)\setminus\overline B(p,r).
\]
For $x\in A_r(p)$, we have $x\in B(0,\frac{1}{2})$,
$q\in B(x,1)$ and $B(x,3)\Subset B(0,8)$.
For a Borel set $A\subset\R^k$, set
$\nu(A)\coloneqq\mu(A\cap B(0,3))$.
The Wolff estimate \eqref{eq:wolff}, with $R=1$, gives
\begin{equation}\label{eq:atom-wolff}
  v(x)\le C v(q)+C\int_0^2
  \nu(B(x,t))^{\frac{1}{k-1}}\,\frac{dt}{t},
\end{equation}
for $C >0$ which depends only on $k$ and uniform bounds of $g$ on $B(0,8)$.
Note that $\nu$ is finite, and its mass agrees with
$\mu$ on every $B(x,t)$ occurring on the right of~\eqref{eq:atom-wolff}.

For $0<t<r$, Tonelli's theorem gives
\begin{align*}
  \fint_{A_r(p)}\nu(B(x,t))\,dx
  &=\frac1{|A_r(p)|}
  \int_{\R^k}|A_r(p)\cap B(y,t)|\,d\nu(y)\\
  &\leq \frac1{|A_r(p)|}
  \int_{B(p,3r)}|B(0,t)|\,d\nu(y)\\
  &=\frac1{2^k-1}(\frac{t}{r})^k\,\nu(B(p,3r)).
\end{align*}
Since $0<\frac{1}{k-1}\le1$, the function $s\mapsto s^{\frac{1}{k-1}}$ is concave on
$[0,\infty)$. Tonelli's theorem and Jensen's inequality give
\begin{align*}
  \fint_{A_r(p)}\int_0^r\nu(B(x,t))^{\frac{1}{k-1}}\,\frac{dt}{t}\,dx
  &\le\int_0^r\left(\fint_{A_r(p)}\nu(B(x,t))\,dx\right)^{\frac{1}{k-1}}\frac{dt}{t}\\
  &\le (2^k-1)^{-\frac{1}{k-1}} \nu(B(p,3r))^{\frac{1}{k-1}}\int_0^r(\frac{t}{r})^{\frac{k}{k-1}}\,\frac{dt}{t}\\
  &=C\nu(B(p,3r))^{\frac{1}{k-1}}.
\end{align*}
For $x\in A_r(p)$ and $r\le t\le2$, the inclusion
$B(x,t)\subset B(p,3t)$ gives
$\nu(B(x,t))^{\frac{1}{k-1}}\le\nu(B(p,3t))^{\frac{1}{k-1}}$.
Averaging \eqref{eq:atom-wolff} over $A_r(p)$ therefore gives
\begin{equation}\label{eq:atom-average}
  \fint_{A_r(p)}v\,dx
  \le C v(q)+C \nu(B(p,3r))^{\frac{1}{k-1}}
  +C \int_r^2\nu(B(p,3t))^{\frac{1}{k-1}}\,\frac{dt}{t}.
\end{equation}
The right side is finite for each $r>0$.

Set $m\coloneqq\nu(\{p\})=\mu(\{p\})$. Then
$\nu(B(p,3t))^{\frac{1}{k-1}}\to m^{\frac{1}{k-1}}$ as $t\to0$.
Fix $\varepsilon>0$ and choose $\delta\in(0,\frac{1}{8})$ such that
$|\nu(B(p,3t))^{\frac{1}{k-1}}-m^{\frac{1}{k-1}}|<\varepsilon$ for $0<t<\delta$.
For $0<r<\delta$, we obtain
\[
  (m^{\frac{1}{k-1}}-\varepsilon)\log(\frac{\delta}{r})
  \le\int_r^\delta\nu(B(p,3t))^{\frac{1}{k-1}}\,\frac{dt}{t}
  \le(m^{\frac{1}{k-1}}+\varepsilon)\log(\frac{\delta}{r}).
\]
Here, $\int_\delta^2\nu(B(p,3t))^{\frac{1}{k-1}}\,\frac{dt}{t}$ is finite and independent of $r$.
Since $\frac{\log(\frac{\delta}{r})}{\log(\frac{1}{r})}\to1$ and $\varepsilon>0$ is arbitrary,
we get
\[
  \lim_{r\to0}\frac1{\log(\frac{1}{r})}
  \int_r^2\nu(B(p,3t))^{\frac{1}{k-1}}\,\frac{dt}{t}=m^{\frac{1}{k-1}}.
\]
The first two terms on the right of \mbox{\eqref{eq:atom-average}},
divided by $\log(\frac{1}{r})$, go to zero as $r\to0$.
Assume also that $2r<r_p$. Then
$A_r(p)\subset B(p,r_p)\setminus\{p\}$, and
assumption~\mbox{\eqref{eq:atom-log-lower}} gives
\[
  \fint_{A_r(p)}v\,dx \ge\alpha\log(\frac{1}{2r})-C_p.
\]
Thus, dividing \mbox{\eqref{eq:atom-average}} by $\log(\frac{1}{r})$ and
taking the limit, we get $\alpha\le Cm^{\frac{1}{k-1}}$, which means
$\mu(\{p\})\ge C^{-(k-1)}\alpha^{k-1}$. This gives the required uniform estimate for
$p \in B(0,\frac{1}{4})$.
A finite cover by corresponding smaller balls gives a uniform constant on every fixed compact
subset of $D$.
\end{proof}

We now use Lemma~\ref{lem:background} to construct a $k$-superharmonic function
from $- \log \lambda$. Then we appeal to the properties of the associated Riesz measure.

\begin{proposition}\label{prop:factor}
Fix $(L^k,g)$ with $k \geq 3$ and let
$\lambda\ge0$ be continuous, smooth on the nonempty set $\mathcal R\coloneqq\{\lambda>0\}$.
Suppose that for every $p\in L$ there exist an open neighborhood
$p\in U\Subset L$ and a constant $C\ge0$ such that
\[
  \Ric_{\lambda^2g}
  \le C\lambda^2g
  \quad\text{on }U\cap\mathcal R.
\]
Then $\mathcal R$ is dense in $L$. If, in addition,
$\lambda\in C^{0,\alpha}_{\mathrm{loc}}(L)$ for some $\alpha\in(0,1]$,
then $Z\coloneqq\{\lambda=0\}$ is discrete.
\end{proposition}
\begin{proof}
Fix $p\in\partial\mathcal R$ and choose $p\in U\Subset L$ as in the hypothesis.
Using Lemma~\ref{lem:background} and the notation of
Remark~\ref{rem:local-conformal-change}, we choose a connected neighborhood
$p\in D\Subset U$ and make a local conformal change of $g$, together with the
corresponding change of $\lambda$, such that
\eqref{eq:ricci-gap} holds on $D\cap\mathcal R$.
This conformal change preserves the condition
$\lambda\in C^{0,\alpha}_{\mathrm{loc}}(D)$.
Set $f\coloneqq\log\lambda$ on $D\cap\mathcal R$. At a point where
$\nabla_g f\ne0$, evaluate \eqref{eq:conformal-ricci} on
$e\coloneqq\frac{\nabla_gf}{|\nabla_gf|_g}$ to get
\[
  Q(e,e)=\Delta_gf+(k-2)\nabla_g^2f(e,e) > 0.
\]
Thus,
\begin{align*}
  \Delta_{k,g}f
  &=\operatorname{div}_g
  \bigl(|\nabla_gf|_g^{k-2}\nabla_gf\bigr)\\
  &=|\nabla_gf|_g^{k-2}\Delta_gf
  +(k-2)|\nabla_gf|_g^{k-4}
  \nabla_g^2f(\nabla_gf,\nabla_gf)\\
  &=|\nabla_gf|_g^{k-2}
  \bigl(\Delta_gf+(k-2)\nabla_g^2f(e,e)\bigr)\\
  &=|\nabla_gf|_g^{k-2}Q(e,e) \geq 0.
\end{align*}
Since $k > 2$, at a point where $\nabla_gf=0$ we also have $\Delta_{k,g}f = 0.$ Thus
$\Delta_{k,g}f\ge0$ everywhere on $D\cap\mathcal R$.

Choose $\Lambda>\sup_D\lambda$, and define the \emph{logarithmic potential}
\[
  v(x)\coloneqq
  \begin{cases}
    \log(\frac{\Lambda}{\lambda(x)}),&x\in D\cap\mathcal R,\\
    +\infty,&x\in D\cap Z.
  \end{cases}
\]
This function is nonnegative and lower semicontinuous.
We now verify the comparison property in the definition of
$k$-superharmonicity on $D$.
Let $V\Subset D$ and let $h\in C(\overline V)$ be $k$-harmonic
on $V$, with $h\le v$ on $\partial V$.
For real $s$, write $s_+\coloneqq\max\{s,0\}$.
For $\varepsilon>0$, set
$w_\varepsilon\coloneqq(h-v-\varepsilon)_+$.
Set $ K_\varepsilon\coloneqq \{q\in\overline V:v(q)\le h(q)-\varepsilon\}.$
Lower semicontinuity of $v-h$ makes $K_\varepsilon$ closed in $\overline V$,
hence compact. Since $h\le v$ on $\partial V$ and $v=+\infty$ on $Z$, we get
$ K_\varepsilon\cap(\partial V\cup Z)=\varnothing,$ and so
$\operatorname{supp}w_\varepsilon\subset K_\varepsilon\subset V\cap\mathcal R.$
Choose an open set $K_\varepsilon\subset\Omega\Subset V\cap\mathcal R$.
Since $h-v-\varepsilon\in W^{1,k}(\Omega)$, we get
\[
  w_\varepsilon\in W^{1,k}_0(\Omega),\qquad
  \nabla_gw_\varepsilon
  =\mathbf1_{\{h-v>\varepsilon\}}(\nabla_gh-\nabla_gv)
  \quad\text{almost everywhere}.
\]

On $\Omega$, $v$ is smooth and satisfies
$\Delta_{k,g}v=-\Delta_{k,g}f\le0$, so $v$ is a smooth supersolution.
Recalling that $h$ is $k$-harmonic and using $w_\varepsilon$ as a test function, we obtain
\begin{align*}
  0&\ge\int_\Omega
  \big\langle |\nabla_gh|_g^{k-2}\nabla_gh-
  |\nabla_gv|_g^{k-2}\nabla_gv,\,
  \nabla_gw_\varepsilon\big\rangle_g\,d\vol_g\\
  &=\int_\Omega\mathbf1_{\{h-v>\varepsilon\}}
  \big\langle |\nabla_gh|_g^{k-2}\nabla_gh-
  |\nabla_gv|_g^{k-2}\nabla_gv,\,
  \nabla_gh-\nabla_gv\big\rangle_g\,d\vol_g
  \ge0,
\end{align*}
where the last inequality is because of strict monotonicity of $\xi\mapsto|\xi|_g^{k-2}\xi$,
applied to $\nabla_g h$ and $\nabla_g v$.
Since the integral now has to vanish, this implies $\nabla_g w_\varepsilon=0$ almost everywhere.
Since $w_\varepsilon$ has compact support, we get $w_\varepsilon\equiv 0$.
Letting $\varepsilon\to0$ proves $h\le v$ on $V$.
Also $v$ is finite somewhere in $D$, since
$p\in\overline{\mathcal R}$. Thus $v$ is $k$-superharmonic on $D$.

Theorem~\ref{thm:potential} now gives the Riesz measure
$\mu\coloneqq-\Delta_{k,g}v$ on $D$ and locally zero $k$-capacity of
$D\cap Z$. In particular $D\cap Z$ has empty interior.
It follows that $\overline{\mathcal R}$ is open, since for every point $p\in\partial\mathcal R$
we just got $D$ contained in $\overline{\mathcal R}$.
It is also closed and nonempty, so connectedness gives
$\overline{\mathcal R}=L$. Hence, every zero of $\lambda$ belongs to $\partial\mathcal R$.
We also showed that every point of $Z \subset \partial\mathcal R$ has a neighborhood $D$ with
$D \cap Z$ of locally zero $k$-capacity, which implies that
$Z$ itself has locally zero $k$-capacity.

Suppose now that $\lambda\in C^{0,\alpha}_{\mathrm{loc}}(L)$, for some $\alpha\in(0,1]$.
Let $p \in Z$. Choosing a connected neighborhood $p \in D$ and performing the conformal
change, as in the beginning of this proof, preserves $\lambda\in C^{0,\alpha}_{\mathrm{loc}}(D)$.
Let $\Lambda,v,\mu$ be defined as before.
Choose $r>0$ so that $B_g(p,2r)\Subset D$ has normal coordinates centered at $p$.
Then for some uniform $C>0$,
\[
  \lambda(x)=|\lambda(x)-\lambda(z)|
  \le C\operatorname{dist}_g(x,z)^\alpha \leq C |x-z|^\alpha
  \qquad
  \bigl(x\in B_g(p,2r),\ z\in Z\cap\overline{B_g(p,r)}\bigr).
\]
Hence, for these $x,z$ with $x\ne z$,
\[
  v(x)\ge\alpha\log\frac1{|x-z|}
  +\log(\frac{\Lambda}{C}).
\]
Lemma~\ref{lem:atom} therefore gives
\[
  \mu(\{z\})\ge m_0\coloneqq c\alpha^{k-1}>0
  \qquad\bigl(z\in Z\cap\overline{B_g(p,r)}\bigr),
\]
with $c$ independent of $z$.
Set $K\coloneqq\overline{B_g(p,r)}$. For any finite
$F\subset K\cap Z$, we have
\[
  m_0\#F\le\sum_{z\in F}\mu(\{z\})\le\mu(K)<\infty.
\]
Thus $K\cap Z$ is finite. Since $p$ lies in the interior of $K$, it is isolated in $Z$.
Since $p\in Z$ was arbitrary, this proves discreteness.
\end{proof}

With this, we are ready to prove Theorem~\ref{thm:A}.
\begin{proof}[Proof of Theorem~\ref{thm:A}]
A nonconstant map has nonempty $\mathcal R_u$. By Lemma~\mbox{\ref{lem:regularity}},
$\lambda \in C^{0,\alpha}_{\mathrm{loc}}$
and $\lambda$ is smooth on $\mathcal R_u$.
Lemma~\mbox{\ref{lem:gauss}} gives the bound, which with Proposition~\mbox{\ref{prop:factor}}
proves that $Z_u$ is discrete.
Note that $Z_u$ is closed because $du$ is continuous.
\end{proof}

\section{Absence of critical points for smooth maps}\label{sec:smooth}
First, we show that under the same local Ricci bound, $\lambda$ cannot
vanish to infinite order at a point on the boundary of $\{\lambda > 0\}$.
This is a standard comparison argument: we apply the maximum
principle on an annulus to compare the maxima of $\lambda$ on
three concentric geodesic spheres.

\begin{lemma}\label{lem:finite-order}
Fix $(L,g)$ and $p\in U\Subset L$. Suppose $\lambda:U\to[0,\infty)$ is
continuous and smooth on $\mathcal R\coloneqq\{x\in U:\lambda(x)>0\}$.
Assume $p\in\{\lambda=0\}\cap\overline{\mathcal R}$ and for some $C\ge0$ we have
\[
  \Ric_{\lambda^2g}\le C\lambda^2g \quad\text{on }U\cap\mathcal R.
\]
Then there are $c,A,r_0>0$ such that
\[
  \max_{\overline B_{g}(p,r)}\lambda\ge cr^A\qquad\text{for $0<r<r_0$}.
\]
\end{lemma}

\begin{proof}
As in the proof of Proposition~\mbox{\ref{prop:factor}}, we choose a connected
neighborhood $p\in D\Subset U$ and make a local conformal change of $g$, together with
the corresponding change of $\lambda$, such that
\mbox{\eqref{eq:ricci-gap}} holds on $D\cap\mathcal R$.
The original and modified metrics are locally uniformly equivalent,
and $\lambda$ also scales by a smooth positive factor,
so it is enough to prove the inequality after this change.

Shrink $D$ to a normal ball centered at $p$. Set $f\coloneqq\log\lambda$ on
$D\cap\mathcal R$. Taking the trace of
\mbox{\eqref{eq:conformal-ricci}} and then applying the chain rule give
\[
  \begin{aligned}
    \Delta_g f+\frac{k-2}{2}|df|_g^2&\ge\frac{kc_0}{2(k-1)},\\
    \Delta_g(\lambda^{\frac{k-2}{2}})
    &\ge\frac{k(k-2)c_0}{4(k-1)}\lambda^{\frac{k-2}{2}}>0
    \quad\text{on }D\cap\mathcal R.
  \end{aligned}
\]
Thus, for every $r>0$ such that $B_g(p,r)\Subset D$, $\lambda^{\frac{k-2}{2}}$ cannot
have a positive local maximum in $B_g(p,r)$.
Since $p\in\overline{\mathcal R}$, we have
$B_g(p,r)\cap\mathcal R\ne\varnothing$.
Hence,
\[
  M(r)\coloneqq\max_{\partial B_g(p,r)}\lambda
  =\max_{\overline B_g(p,r)}\lambda>0,\qquad M(0):=0.
\]
Clearly $M(r)$ is nondecreasing. We quickly show that $M(r)$ is continuous on $[0,\rho)$
for $\rho>0$ with $\overline B_g(p,\rho)\Subset D$. Fix $\varepsilon>0$.
Uniform continuity of $\lambda$ on $\overline B_g(p,\rho)$
gives a $\delta>0$ such that
$|\lambda(x)-\lambda(y)|<\varepsilon$ whenever
$x,y\in\overline B_g(p,\rho)$ and $\operatorname{dist}_g(x,y)<\delta$.
For $0<r,s<\rho$ with $|r-s|<\delta$, we have
$\operatorname{dist}_g(\exp_p(r\theta),\exp_p(s\theta)) = |r-s|<\delta$
for every unit vector $\theta\in T_pL$. Consequently,
$|M(r)-M(s)|\le\max_{|\theta|_{g_p}=1}
|\lambda(\exp_p(r\theta))-\lambda(\exp_p(s\theta))|<\varepsilon$.
Combining this with $M(r)\to0 = M(0)$ as $r\to0$, we get the required continuity.

On $B_g(p,\rho)$, define $r(x)\coloneqq\operatorname{dist}_g(p,x)$.
On $B_g(p,\rho)\setminus\{p\}$, Petersen~\cite[Exercise~5.9.27(3)]{Petersen} gives
\[
  \nabla_g^2r=\frac1r(g-dr\otimes dr)+O(r).
\]
Taking the trace and using $|dr|_g=1$, we obtain
\[
  \Delta_g r=\frac{k-1}{r}+O(r).
\]
Thus there is a constant $C\ge0$ such that
\begin{equation}\label{eq:distance-laplacian}
  \frac{\Delta_g r}{k-1}\le\frac1r+Cr
  \qquad(0<r<\rho).
\end{equation}
For $0<r\le\rho$, define
\begin{equation}\label{eq:corrected-radius-definition}
  t(r)\coloneqq\log r+
  \int_0^r\frac{e^{-\frac{Cs^2}{2}}-1}{s}\,ds.
\end{equation}
The integrand in~\eqref{eq:corrected-radius-definition} is $O(s)$ as $s\to0$, so
$t(r)=\log r+O(r^2)$.
Moreover,
\begin{equation}\label{eq:corrected-radius}
  t'(r)=\frac{e^{-\frac{Cr^2}{2}}}r>0,\qquad
  t''(r)+(r^{-1}+Cr)t'(r)=0.
\end{equation}

Set $F(r)\coloneqq\log M(r)$ for $0<r<\rho$. Fix $0<r_1<r_2<\rho$
and let $\ell:\mathbb R\longrightarrow\mathbb R$ be the affine function defined by
$\ell(t(r_i))=F(r_i)$
for $i=1,2$. The slope of $\ell$ is $a\coloneqq\frac{F(r_2)-F(r_1)}{t(r_2)-t(r_1)}\ge0$, since
$M$ is nondecreasing and $t$ is strictly increasing.
We claim that $F(r)\le\ell(t(r))$ for $r_1\le r\le r_2$.
On $\mathcal A\coloneqq\{x\in B_g(p,\rho):r_1\le r(x)\le r_2\}$,
the function $x\mapsto\lambda(x)\exp(-\ell(t(r(x))))$ is continuous
and at most $1$ on $\partial B_g(p,r_1)$ and $\partial B_g(p,r_2)$.
If the claim fails, this function
attains a maximum greater than $1$ at some $q\in\mathcal A$ with $r_1<r(q)<r_2$.
Then $\lambda(q)>0$, so $f-\ell(t(r))$ is smooth near $q$ and has a
local maximum at $q$.

Set $e\coloneqq\nabla_g r(q)$. Since $|\nabla_g r|_g=1$ away from $p$, we have
$|e|_g=1$ and $\nabla_g^2r(e,e)=0$.
Since $f-\ell(t(r))$ has a local maximum at $q$, we have
\[
  \nabla_g f=at'e,\qquad
  \nabla_g^2f(e,e)\le at'',\qquad
  \Delta_g f\le a(t''+t'\Delta_g r),
\]
where $t'$ and $t''$ denote the derivatives at $r(q)$.
Thus $df(e)^2=a^2(t')^2=|df|_g^2$. Evaluating \eqref{eq:ricci-gap}
and \eqref{eq:conformal-ricci} on $(e,e)$ at $q$, we obtain
\[
  \begin{aligned}
    c_0\le Q(e,e)
    &=(k-2)\nabla_g^2 f(e,e)+\Delta_g f\\
    &\le a\bigl((k-1)t''+t'\Delta_g r\bigr)\\
    &\le a(k-1)\bigl(t''+(r(q)^{-1}+Cr(q))t'\bigr)
    &&\text{(by \eqref{eq:distance-laplacian}, since $at'\ge0$)}\\
    &=0
    &&\text{(by \eqref{eq:corrected-radius})}.
  \end{aligned}
\]
This contradicts $c_0>0$ and proves the claim.

Since $r_1,r_2$ were arbitrary, the claim says that
$G\coloneqq F\circ t^{-1}$ is convex on $(-\infty,t(\rho))$.
As $r\to0$, we have $t(r)\to-\infty$ and $M(r)\to0$, so
$G(\tau)\to-\infty$ as $\tau\to-\infty$. We may therefore choose
$\tau_1<\tau_2<t(\rho)$ with $G(\tau_1)<G(\tau_2)$. Convexity gives
\[
  \frac{G(\tau_1)-G(\tau)}{\tau_1-\tau}
  \le A\coloneqq\frac{G(\tau_2)-G(\tau_1)}{\tau_2-\tau_1}
  \qquad(\tau<\tau_1).
\]
Thus $A>0$ and $\log M(r)\ge G(\tau_1)+A(t(r)-\tau_1)$ whenever
$t(r)<\tau_1$. Since $t(r)-\log r=O(r^2)$, choose
$0<r_0<t^{-1}(\tau_1)$ such that $|t(r)-\log r|\le1$ for every
$0<r<r_0$. Then $t(r)<\tau_1$ and
\[
  \log M(r)\ge G(\tau_1)+A(\log r-1-\tau_1)
  \qquad(0<r<r_0).
\]
Exponentiating gives
\[
  M(r)\ge cr^A\quad(0<r<r_0),
  \qquad c\coloneqq e^{G(\tau_1)-A(1+\tau_1)}>0,
\]
which is what we needed to prove.
\end{proof}

Next, we use~\cite[Theorem~B]{Branch}, which gives the following estimate.

\begin{lemma}\label{lem:taylor}
Let $u:(L^k,g)\to(M^n,h)$ be smooth and weakly-conformal, with
$n\ge k\ge3$. For every $p\in\Crit(u)$ and every integer $j\ge1$,
there exist constants $C,r>0$ such that
\[
  |du(x)|_{g,h}\le C\operatorname{dist}_g(p,x)^j
  \qquad(x\in B_g(p,r)).
\]
\end{lemma}

With this, we are ready to prove Theorem~\ref{thm:B}.
\begin{proof}[Proof of Theorem~\ref{thm:B}]
If $\mathcal R_u = \emptyset$, then $du=0$ everywhere and so $u$ is constant. Otherwise, if
$\mathcal R_u\ne L$, choose
$p\in\partial\mathcal R_u$. Such a point exists because $L$ is
connected.
Lemma~\mbox{\ref{lem:gauss}} gives a uniform Ricci bound near $p$.
Thus, Lemma~\mbox{\ref{lem:finite-order}} gives $c,A,r_0>0$ such that
$\max_{\overline B_g(p,r)}\lambda\ge cr^A$ for $0<r<r_0$.
On the other hand, Lemma~\mbox{\ref{lem:taylor}} gives
$|du(x)|_{g,h}=O(\operatorname{dist}_g(p,x)^j)$ for every integer $j\ge1$.
Thus $\max_{\overline B_g(p,r)}\lambda=O(r^j)$ for every integer $j\ge1$.
Choosing $j>A$ gives a contradiction. Hence $\mathcal R_u=L$,
and weak conformality gives that $u$ is an immersion.
\end{proof}

\section{Two-map unique continuation}\label{sec:continuation}

For each class, we first write the equation satisfied by the mean curvature of $u|_{\mathcal R_u}$.
For every $y\in M$ and every oriented $k$-plane $E\subset T_yM$, let $\tau_E$ be its
oriented unit simple $k$-vector. In class~(i), define the normal
vector $\mathcal H(y,E)\in E^\perp$ by
\[
  h(\mathcal H(y,E),X)\coloneqq-\frac1k\,d\phi_y(X,\tau_E),
  \qquad X\in E^\perp.
\]
In class~(ii), set $\mathcal H\coloneqq0$.
Then $\mathcal H(y,E)\in E^{\perp}\subset T_yM$ depends smoothly on \(y\) and \(E\subset T_yM\).
If $V\subset\mathcal R_u$ is open and $u|_V$ is an embedding, then
$\Sigma\coloneqq u(V)$ has averaged mean curvature
$H(y)=\mathcal H(y,T_y\Sigma)$ for $y\in\Sigma$, by Lemma~\ref{lem:regularity}.
In class~(ii), $\mathcal H(y,E)=0$ for either orientation of $E$.

We write the standard mean-curvature equation for a graph and estimate
the difference between two solutions.
Next, we apply Aronszajn's unique continuation theorem~\cite[Remark~3]{Aronszajn}.

\begin{lemma}\label{lem:graph}
Let $W\subset\R^k\times\R^{n-k}$ be the image of a target chart,
and let $h$ be a smooth Riemannian metric on $W$.
Fix $x_0\in\R^k$ and $r>0$, and set $D\coloneqq B(x_0,r)$.
Let $f\in C^\infty(D,\R^{n-k})$ satisfy $(x,f(x))\in W$ for every $x\in D$.
Define
\[
  \begin{gathered}
    G_f:D\longrightarrow W,\qquad G_f(x)\coloneqq(x,f(x)),\\
    E_f(x)\coloneqq dG_f|_x(\R^k),\qquad q_f\coloneqq G_f^*h.
  \end{gathered}
\]
Fix an orientation of $D$ and give $E_f(x)$ the orientation induced by $dG_f|_x$.
Write $H_f(x)\in E_f(x)^{\perp_h}$ for the averaged mean curvature
of $G_f(D)$ at $G_f(x)$.
If $H_f(x)=\mathcal H(G_f(x),E_f(x))$ for every $x\in D$, then
\begin{equation}\label{eq:graph}
  (q_f)^{ij}(x)\partial_i\partial_j f^\alpha(x)
  =B^\alpha(x,f(x),Df(x))
  \qquad(x\in D,\ 1\le\alpha\le n-k),
\end{equation}
where $B:W\times\R^{(n-k)\times k}\to\R^{n-k}$ is smooth and
depends only on $h$, $\mathcal H$ and the chosen orientation of $D$.
Suppose $f_1,f_2$ solve \eqref{eq:graph} on $D$ with the same $h$
and the same $B$. Set
$w\coloneqq f_1-f_2$ and $a^{ij}\coloneqq (q_{f_1})^{ij}$.
For every open ball $D'\Subset D$, there is a constant $C>0$ such that on $D'$
\begin{equation}\label{eq:aronszajn}
  |a^{ij}\partial_i\partial_jw|
  \le C\bigl(|w|+|Dw|\bigr).
\end{equation}
Here, $(a^{ij})$ is smooth, symmetric and uniformly positive definite on $D'$.

If $k=n$, then we have $\mathbb{R}^{n-k}=\{0\}$, so $f\equiv0$ and $G_f(D)=D$.
\end{lemma}
\begin{proof}
Let $\Gamma^A_{BC}$ be the Christoffel symbols of $h$ on $W$.
For $z\in W$, define $\Gamma_h(z):\R^n\times\R^n\to\R^n$ by
\[
  \Gamma_h(z)(v_1,v_2)^A\coloneqq\Gamma^A_{BC}(z)v_1^Bv_2^C
  \qquad(v_1,v_2\in\R^n,\ 1\le A\le n).
\]
We work on $D$. Let $P_f^\perp:\R^n\to E_f^{\perp_h}$ be the $h_{G_f}$-orthogonal projection.
Define
\[
  T_f:\R^{n-k}\longrightarrow E_f^{\perp_h},\qquad
  T_fz\coloneqq P_f^\perp(0,z).
\]
If $T_fz=0$, then $(0,z)\in E_f$, so
$(0,z)=(v,df(v))$ for some $v\in\R^k$. Hence $v=0$ and $z=0$.
Both $\R^{n-k}$ and $E_f^{\perp_h}$ have the same dimension, so $T_f$ is an isomorphism.
Note that $T_f$ and $T_f^{-1}$ depend smoothly on $(x,f(x),Df(x))$.

By definition,
\[
  kH_f=P_f^\perp\!\left[
  (q_f)^{ij}\left(\partial_i\partial_jG_f+
  \Gamma_h(G_f)(\partial_iG_f,\partial_jG_f)\right)\right].
\]
Since $\partial_i\partial_jG_f=(0,\partial_i\partial_jf)$,
applying $T_f^{-1}$ gives \mbox{\eqref{eq:graph}}, with
\[
  \begin{aligned}
    B(x,f(x),Df(x))\coloneqq T_f^{-1}\!\big[&
    k\mathcal H(G_f,E_f) -P_f^\perp\big((q_f)^{ij}
    \Gamma_h(G_f)(\partial_iG_f,\partial_jG_f)\big)\big].
  \end{aligned}
\]
The right-hand side depends only on $(x,f(x),Df(x))$ and is smooth in them.
Replacing $f(x)$ and $Df(x)$ by independent
variables $y$ and $P$ defines the smooth map $B(x,y,P)$, where
$(x,y)\in W$ and $P\in\R^{(n-k)\times k}$.

Subtracting \mbox{\eqref{eq:graph}} for $f_1$ and $f_2$ gives on $D$
\begin{equation}\label{eq:graph-difference}
  \begin{aligned}
    a^{ij}\partial_i\partial_jw
    ={}&B(x,f_1(x),Df_1(x))-B(x,f_2(x),Df_2(x))\\
    &+((q_{f_2})^{ij}-(q_{f_1})^{ij})\partial_i\partial_jf_2.
  \end{aligned}
\end{equation}
Fix an open ball $D'\Subset D$. The set
\[
  K\coloneqq\bigcup_{\ell=1}^2
  \{(x,f_\ell(x),Df_\ell(x)):x\in\overline{D'}\}
\]
is a compact subset of $W\times\R^{(n-k)\times k}$.
Replacing $f(x)$ and $Df(x)$ by $y$ and $P$ in the formula for
$(q_f(x))^{-1}$ defines a smooth matrix-valued function of $(x,y,P)$.
Its smoothness and that of $B$ give some $C>0$ such that on $D'$
\begin{equation}\label{eq:graph-coefficient-estimate}
  \begin{aligned}
    &|B(x,f_1(x),Df_1(x))-B(x,f_2(x),Df_2(x))|\\
    &\qquad{}+|((q_{f_2})^{ij}-(q_{f_1})^{ij})_{i,j}|
    \le C\bigl(|w|+|Dw|\bigr).
  \end{aligned}
\end{equation}
Combining \eqref{eq:graph-difference} and
\eqref{eq:graph-coefficient-estimate} with the bound
$\sup_{\overline{D'}}|D^2f_2|<\infty$ proves
\eqref{eq:aronszajn}.
Since $(a^{ij})=((q_{f_1})_{ij})^{-1}$ and $q_{f_1}$ is a smooth metric, compactness of
$\overline{D'}$ gives a constant
$c>0$ such that on $D'$, for every $\xi\in\R^k$,
\[
  c^{-1}|\xi|^2\le a^{ij}\xi_i\xi_j\le c|\xi|^2.
\]
\end{proof}

Once we get that the local images of $u$ and $v$ agree, we write $v=u\circ F$ for a
conformal diffeomorphism $F$ between neighborhoods of $p$ in $L$.
We use the following lemma to show that $F$ is the identity near $p$.
For a proof of the lemma,
see~\cite[Chapter~IV, Section~6, and Chapter~I,
proof of Theorem~3.2]{KobayashiTG}.

\begin{lemma}\label{lem:two-jet}
Let $F$ be a smooth conformal local diffeomorphism between open
subsets of $(L^k,g)$, $k\ge3$. If its two-jet at $p$ equals that
of the identity, then $F$ is the identity near $p$.
\end{lemma}

We now prove that two maps agree near a common regular point if their
values and all their derivatives agree there.

\begin{lemma}\label{lem:regular-uc}
Let $u,v$ belong to the same class (i) or (ii) with the same background data,
and let $p\in\mathcal R_u\cap\mathcal R_v$.
Suppose $u(p)=v(p)$ and, in common local coordinates at $p$ and $u(p)$,
\begin{equation}\label{eq:equal-jets}
  \partial^\beta u(p)=\partial^\beta v(p),
  \qquad\text{for every $\beta$ with $|\beta|\ge1$.}
\end{equation}
Then $u=v$ on a neighborhood of $p$.
\end{lemma}
\begin{proof}
By Lemma~\mbox{\ref{lem:regularity}}, $u$ and $v$ are smooth near $p$. Note that by the chain rule,
the condition~\eqref{eq:equal-jets} is independent of local coordinates.
We choose target local coordinates centered at $u(p)$ in which
$du_p(T_pL)=\R^k\times\{0\}$. For $s\in\{u,v\}$, write $s=(X_s,Y_s)$ in $\R^k\times\R^{n-k}$.
Since $d(X_u)_p=d(X_v)_p$ is invertible, each $X_s$ maps a
neighborhood of $p$ diffeomorphically onto the same open ball $D\subset\mathbb R^k$
centered at $0$. Define $f_s\coloneqq Y_s\circ X_s^{-1}$ on $D$.
Then, for every $\beta$ with $|\beta| \geq 1$ we have:
\begin{align*}
\partial^\beta(f_u\circ X_v)(p)
&=\partial^\beta(f_u\circ X_u)(p)
=\partial^\beta Y_u(p)\\
&=\partial^\beta Y_v(p)
=\partial^\beta(f_v\circ X_v)(p).
\end{align*}
Composing with $X_v^{-1}$ gives
$\partial^\beta(f_u-f_v)(0)=0$.
In class~(i), $X_u$ and $X_v$ induce the same orientation on $D$
because $d(X_u)_p=d(X_v)_p$ and $D$ is connected.
Fix an open ball $D'\Subset D$ centered at $0$. If $n>k$,
Lemma~\mbox{\ref{lem:graph}} gives \mbox{\eqref{eq:aronszajn}} on $D'$ for
$w\coloneqq f_u-f_v$. Aronszajn's unique continuation
theorem~\mbox{\cite[Remark~3]{Aronszajn}} gives $w=0$ on $D'$.
If $n=k$, then $f_u=f_v=0$.

Define $F\coloneqq X_u^{-1}\circ X_v$ on $X_v^{-1}(D')$.
Since $f_u=f_v$ on $D'$, we have $v=u\circ F$ on $X_v^{-1}(D')$. The equations
$u^*h=\lambda_u^2g$ and $v^*h=\lambda_v^2g$ give
\[
  F^*g=\left(\frac{\lambda_v}{\lambda_u\circ F}\right)^2g.
\]
Thus $F$ is conformal. Since $X_v$ and $X_u$ have the same two-jet
at $p$, so do $F$ and $X_u^{-1}\circ X_u=\operatorname{id}$.
Lemma~\ref{lem:two-jet} gives $F=\operatorname{id}$ on a neighborhood
of $p$, so there $v=u\circ F=u$.
\end{proof}

We now prove Theorem~\ref{thm:C}.

\begin{proof}[Proof of Theorem~\ref{thm:C}]
If one map is constant and the other is not, agreement on an open
set contradicts Theorem~\ref{thm:A} for the nonconstant map, since its
differential would vanish on that set. If both maps are constant, then if they
agree just at one point, they agree everywhere.

Now assume that both maps are nonconstant. By
Theorem~\mbox{\ref{thm:A}}, the set
$Z\coloneqq Z_u\cup Z_v$ is closed and discrete.
Since $L$ is connected and $k\ge3$,
$R\coloneqq L\setminus Z$ is path connected and dense in $L$.
Define
\[
  E\coloneqq\{p\in R:\ u=v
  \text{ on some open neighborhood of }p\}.
\]
It is open in $R$. Moreover, $E \neq \emptyset$ because the open set from the
assumption of the theorem on which they agree intersects the dense set $R$.
It remains to show that $E$ is closed in $R$.

Take a sequence $p_j\in E$ which converges to some $p\in R$.
By continuity, we have $u(p)=v(p)$. Both $u$ and $v$ are smooth on $R$ and agree
near each $p_j$, so all their partial derivatives agree at $p$ by
continuity. Hence by Lemma~\mbox{\ref{lem:regular-uc}}, they agree on a neighborhood
of $p$, and thus $p\in E$. Since $R$ is connected, we conclude that $E=R$.
Since $u$ and $v$ are continuous and $R$ is dense in $L$, we conclude that $u=v$ on all of $L$.
\end{proof}

\section{Discreteness of fibres}\label{sec:fibres}

We make the following definition and then prove Theorem~\ref{thm:D}.
For each open set $U\Subset L$, define the finite Borel measure
$\mu_U$ on $M$ by
\[
  \mu_U(A)\coloneqq
  \int_{U\cap u^{-1}(A)}\lambda^k\,d\vol_g,
  \qquad\text{for every Borel set $A\subset M$.}
\]

\begin{proof}[Proof of Theorem~\ref{thm:D}]
Fix $y\in M$ and $p\in u^{-1}(y)$. If $p$ is regular, $u$ is an
immersion near $p$, so $p$ is isolated in its fibre.

Suppose $p$ is critical. By Theorem~\ref{thm:A}, choose a normal ball
$B_g(p,R)\Subset L$ in which $p$ is the only critical point.
For $0<a<b<R$, the set
$\{q\in u^{-1}(y):a\le \operatorname{dist}_g(p,q)\le b\}$ is compact
and consists of isolated points, so it is finite. Thus only finitely many
spheres of radii in $[a,b]$ about $p$ intersect the fibre. We can therefore
choose $0<r<R$ and $0<r_{j+1}<r_j<r$ with $r_j\to0$ such that
$y\notin u(\partial B_g(p,r))\cup u(\partial B_g(p,r_j))$ for every $j$.
Set $U\coloneqq B_g(p,r)$ and $A_j\coloneqq U\setminus\overline{B_g(p,r_j)}$.

Choose $s_0>0$ so that $B_h(y,2s_0)\Subset M$ is a normal ball
with $B_h(y,2s_0) \cap u(\partial U) = \emptyset$. On $\overline{A_j}$, $u$
is a smooth immersion. Define
\[
  b_j\coloneqq\operatorname{Vol}_{u^*h}\bigl(\partial B_g(p,r_j)\bigr)
  =\int_{\partial B_g(p,r_j)}\lambda^{k-1}\,dA_g,
\]
where $dA_g$ is the measure induced by $g$.
Since $\lambda$ is bounded on $\overline U$ and
$\operatorname{Vol}_g(\partial B_g(p,r_j))=O(r_j^{k-1})$,
we have $b_j=O(r_j^{k-1})\to0$.

For $0<t<s_0$, set
$Q_j(t)\coloneqq t^{-k}\mu_{A_j}(\overline{B_h(y,t)})$.
Apply the monotonicity inequality
\cite[Example~2.4, Lemmas~2.7 and~3.3, and (3.7)]{Scharrer}
to $u|_{\overline{A_j}}$. Lemma~\mbox{\ref{lem:regularity}} makes $|H|_h$ bounded
on $U\setminus\{p\}$, and $\sec_h$ is bounded on $B_h(y,2s_0)$. Thus, after decreasing $s_0$ if
necessary, there exists $c\ge0$, independent of $j$, such that
\begin{equation}\label{eq:annular-monotonicity}
  Q_j(s)\le Q_j(t)+c\int_s^t Q_j(\rho)\,d\rho
  +b_j\int_s^t\rho^{-k}\,d\rho
  \qquad(0<s<t<s_0).
\end{equation}
For $0<\rho<s_0$, set $Q(\rho)\coloneqq\rho^{-k}\mu_U(\overline{B_h(y,\rho)})$.
Since $A_j\uparrow U\setminus\{p\}$ and $\lambda(p)=0$, we have
$Q_j(\rho)\uparrow Q(\rho)$. For fixed $0<s<t<s_0$, the last term
in~\eqref{eq:annular-monotonicity}
tends to zero as $j\to\infty$, so monotone convergence gives
$Q(s)\le Q(t)+c\int_s^t Q(\rho)\,d\rho$. Gr\"onwall's inequality gives
$Q(s)\le e^{c(t-s)}Q(t)$. Taking $t\coloneqq \frac{s_0}{2}$ and
$C\coloneqq e^{ct}Q(t)$ yields $\mu_U(B_h(y,s))\le Cs^k$ for $0<s<t$.
Let $\omega_k\coloneqq |B_{\R^k}(0,1)|$.
If $y\notin u(\overline{A_j})$, compactness gives
$B_h(y,s)\cap u(\overline{A_j})=\varnothing$ for all sufficiently small $s>0$,
so $\#(A_j\cap u^{-1}(y))=0$. In the other case, since $\overline{A_j}$ is compact,
we apply the density formula~\cite[Lemma~2.3, identity~(2.3)]{Scharrer}
to $u|_{\overline{A_j}}$. Since $y\notin u(\partial A_j)$, this gives
\[
  \#\bigl(A_j\cap u^{-1}(y)\bigr)
  =\lim_{s\to0}\frac{\mu_{A_j}(B_h(y,s))}{\omega_k s^k}
  \le \frac{C}{\omega_k}.
\]
As $A_j$ increases to $U\setminus\{p\}$, this inequality makes
$U\cap u^{-1}(y)$ finite. Thus $p$ is isolated in its fibre.
Every fibre is closed by continuity of $u$.
\end{proof}

\section{Density of injective points}\label{sec:injective}

McDuff--Salamon~\mbox{\cite[Proposition~2.5.1]{McDuffSalamon}} prove density of
injective points for simple $J$-holomorphic curves. For a minimal immersion
of a closed connected manifold, White~\mbox{\cite[after Theorem~1]{White}}
observes that one injective point is enough to conclude their density. Here, we use
Lemma~\mbox{\ref{lem:graph}} again to compare local graphs of $u$
away from critical points.

\begin{proof}[Proof of Theorem~\ref{thm:E}]
Set $F\coloneqq u^{-1}(u(Z_u))$ and $R\coloneqq L\setminus F$.
By Theorem~\mbox{\ref{thm:A}}, $Z_u$ is closed and discrete, and hence finite by compactness of $L$.
By Theorem~\mbox{\ref{thm:D}} each fibre is finite. Thus, $F$ is finite.
Since $L$ is connected and $k\ge3$, $R$ is connected and dense in $L$.
If $p\in\Inj(u)$, then $p \in \mathcal R_u$ and $u^{-1}(u(p))=\{p\}$.
Hence $u(p)\notin u(Z_u)$, so $p\in R$.

Fix $p\in R$ and write $y\coloneqq u(p)$. Let
$u^{-1}(y) = \{p_1\coloneqq p,p_2,\ldots,p_m\} \subset R$.
Note that $u$ is a smooth immersion on $R$. We want to write the image of $u$ near each $p_i$
as a graph over the same ball in $\R^k$. For each $i$, since $du_{p_i}$
is injective, the linear maps $A:T_yM\to\R^k$ for which $A\circ du_{p_i}$
is invertible form an open dense set. Choose $A$ in the intersection of
these finitely many open dense sets.
Write $A=(A^1,\ldots,A^k)$ and extend them to a basis of $T_y^*M$.
Choose target local coordinates $z=(z^1,\ldots,z^n)$ on a neighborhood $W$ of $y$
such that $z(y)=0$ and $d(z^1,\ldots,z^k)_y=A$. Write
\[
  (X,Y)\coloneqq z\circ u:u^{-1}(W)\longrightarrow\R^k\times\R^{n-k}.
\]
Then $X=(z^1,\ldots,z^k)\circ u$, so by the chain rule,
$dX_{p_i}=A\circ du_{p_i}$, which is invertible for every $i$.

By the inverse function theorem, we can choose disjoint open neighborhoods
$W_i\subset R$ of $p_i$ such that $X|_{W_i}$ is a diffeomorphism onto
an open neighborhood of $0$ in $\R^k$. Choose a ball $D$ centered at $0$
with $D\subset\bigcap_{i=1}^m X(W_i)$ and set $U_i\coloneqq W_i\cap X^{-1}(D)$.
Thus, $u(U_i)$ is the graph of
\[
  f_i\coloneqq Y\circ(X|_{U_i})^{-1}:D\to\R^{n-k}.
\]
By compactness of $L$, the set
\[
  V\coloneqq M\setminus u\Bigl(L\setminus\bigcup_{i=1}^m U_i\Bigr)
\]
is an open neighborhood of $y$ and $u^{-1}(V)\subset\bigcup_i U_i$.
When $p\in\Inj(u)$, we have $m=1$. Then $u|_{U_1}$ is an embedding
and $u^{-1}(V)\subset\Inj(u)\cap U_1$. This proves that $\Inj(u)$ is open and
the claim about neighbourhoods from the statement of the theorem.

Now, we show that $\Inj(u)$ is dense. Assume that $p\in\overline{\Inj(u)}$.
For $2\le i\le m$, let
\[
  E_i\coloneqq\{x\in D:f_i(x)=f_1(x)\}.
\]
If $E_i$ contains a nonempty open set, then graphs of $f_1$ and $f_i$ agree on this open set.
In class~(i), $\phi$ determines the orientation on their common tangent
planes, so both graphs induce the same orientation on $D$. Hence, we get
\[
  H_{f_j}(x)=\mathcal H(G_{f_j}(x),E_{f_j}(x))
  \qquad (x\in D,\ j\in\{1,i\}),
\]
where $d\phi$ determines $\mathcal H$ as in Section~\mbox{\ref{sec:continuation}}.
In class~(ii), $H_{f_1}=H_{f_i}=0$. Thus, in both cases, $f_1$ and $f_i$
satisfy the same equation~\mbox{\eqref{eq:graph}} on $D$. For $n>k$, Lemma~\mbox{\ref{lem:graph}}
and Aronszajn's unique continuation theorem give $f_i=f_1$ on $D$.
For $n=k$, this equality holds because $\R^{n-k}=\{0\}$.
If $f_i=f_1$ on $D$, every $q\in U_1$ has $u(q)=u(q')$ for some
$q'\in U_i$. This contradicts $p\in\overline{\Inj(u)}$.
Thus each $E_i$ is closed with empty interior.

For $q\in U_1\cap u^{-1}(V)$, every preimage of $u(q)$ lies in one of
the $U_i$. Hence $q\in\Inj(u)$ exactly when
$X(q)\notin\bigcup_{i=2}^m E_i$. This finite union has dense
complement in $D$. Since $X|_{U_1}$ is a diffeomorphism onto $D$,
injective points are dense in $U_1\cap u^{-1}(V)$, which is a neighborhood of $p$.
Since $\varnothing\ne\Inj(u)\subset R\cap\overline{\Inj(u)}$, this set
is nonempty. Being open and closed in the connected $R$, it equals $R$.
Density of $R$ in $L$ finishes the proof.
\end{proof}

\end{document}